\documentclass[12pt,twoside]{article}
\usepackage{amsmath}
\usepackage[all]{xy}
\newcommand{\VolumeHeader}{}
\newcommand{\VolumeSerial}{LNS}

\newcommand{\ActivityName}{ {\normalsize {\it 
Summer School on High-dimensional Manifold Topology }}}
\newcommand{\ActivityDate}{ {\normalsize {\it
Trieste, 21 May -- 8 June 2001
}}}

\usepackage{amssymb}
\newcommand{\be}{\begin{equation}}
\newcommand{\ee}{\end{equation}}
\newcommand{\bea}{\begin{eqnarray}}
\newcommand{\eea}{\end{eqnarray}}

\def\Q{\mathbb Q}
\def\R{\mathbb R}
\def\Z{\mathbb Z}
\def\qed{$\square$}
\newcommand{\LectureHeader}{Algebraic surgery}

\begin{document}
\pagestyle{myheadings}
\markboth{\LectureHeader}{\VolumeHeader}
\markright{\VolumeHeader}


\begin{titlepage}


\title{Foundations of algebraic surgery}

\author{Andrew Ranicki\thanks{aar@maths.ed.ac.uk}
\\[1cm]
{\normalsize
{\it Department of Mathematics and Statistics}}\\
{\normalsize
{\it University of Edinburgh, Scotland, UK} }
\\[10cm]
{\normalsize {\it Lecture given at the: }}
\\
\ActivityName 
\\
\ActivityDate 
\\[1cm]
{\small \VolumeSerial} 
}
\date{}
\maketitle
\thispagestyle{empty}
\end{titlepage}

\baselineskip=14pt
\newpage
\thispagestyle{empty}


\begin{abstract}
The algebraic theory of surgery on chain
complexes $C$ with Poincar\'e duality 
$$H^*(C)~ \cong~ H_{n-*}(C)$$
describes geometric surgeries on the chain level. 
The algebraic effect of a geometric surgery on an $n$-dimensional
manifold $M$ is an algebraic surgery on the $n$-dimensional symmetric
Poincar\'e complex $(C,\phi)$ over $\Z[\pi_1(M)]$ with the homology 
of the universal cover $\widetilde{M}$
$$H_*(C)~=~H_*(\widetilde{M})~.$$
The algebraic effect of a geometric surgery on an $n$-dimensional
normal map $(f,b):M\to X$ is an algebraic surgery on the kernel
$n$-dimensional quadratic Poincar\'e complex $(C,\psi)$ over $\Z[\pi_1(X)]$ 
with homology 
$$H_*(C)~=~K_*(M)~=~{\rm ker}(f_*:H_*(\widetilde{M}) \to
H_*(\widetilde{X}))~.$$
For $n>4$ and $i$-connected $(f,b)$ with $2i \leqslant n$ there is a one-one
correspondence between geometric surgeries on $(f,b)$ killing elements
$x \in K_i(M)$ and algebraic surgeries on $(C,\psi)$ killing $x \in
H_i(C)$.  The Wall surgery obstruction of an $n$-dimensional normal map
$(f,b):M \to X$ 
$$\sigma_*(f,b) \in L_n(\Z[\pi_1(X)])$$ 
was originally defined by first making $(f,b)$ $[n/2]$-connected by
geometric surgery below the middle dimension, using forms for even $n$
and automorphisms of forms for odd $n$.
The algebraic theory of surgery identifies $\sigma_*(f,b)$ with the 
cobordism class of the kernel quadratic Poincar\'e complex $(C,\psi)$,
so the algebraic surgery obstruction has the same formulation for
odd and even $n$.
The identification is used for $n=2i$ (resp.  $2i+1$)
to find a representative form (resp. automorphism) without preliminary
geometric surgeries.

\end{abstract}

\vspace{6cm}


{\it AMS numbers:} 57R67, 55U15, 19J25

\newpage
\thispagestyle{empty}
\tableofcontents

\newpage
\setcounter{page}{1}

\section{Introduction}

We compare the homology and chain level descriptions of surgery 
on a manifold, using a minimum of algebraic development.

Manifolds $M$ are to be finite-dimensional, compact, and oriented
(unless stated otherwise), with $C(M)$ denoting the cellular chain complex
for some $CW$ structure on $M$.

Cobordisms $(W;M,M')$ are to be oriented (unless stated otherwise) with
$\partial W=M \cup -M'$, where $-M'$ denotes $M'$ with the opposite
orientation.

\subsection{Background}

The surgery method of classifying manifolds within a homotopy type was
first applied by Kervaire and Milnor \cite{KervaireMilnor} to exotic
spheres, using exact sequences to describe the homology effect of
geometric surgery.  Homology was quite adequate for the subsequent
development of surgery theory on simply-connected manifolds (Browder
\cite{Browder}, Novikov).  Wall \cite{Wall} used a combination of
topology and homology to describe the effect of surgery 
on non-simply-connected manifolds. In general, the homology 
$\Z[\pi_1(M)]$-modules $H_*(\widetilde{M})$ of the universal cover 
$\widetilde{M}$ of a compact manifold $M$ are not finitely generated,
so a chain level approach is indicated.
The algebraic theory of surgery of Ranicki 
\cite{Ranicki1980I},\cite{Ranicki1980II} provided
a model for surgery using chain complexes with Poincar\'e duality.

Surgery was originally developed for differentiable manifolds, but has
since been extended to $PL$ and topological manifolds. The algebraic
theory of surgery applies to all categories of manifolds.

\subsection{The algebraic effect of a geometric surgery}\label{vs}

Let $M$ be an $n$-dimensional manifold. 
Surgery on $S^i \times D^{n-i} \subset M$ results in an $n$-dimensional 
manifold
$$M'~=~(M \backslash S^i \times D^{n-i})\cup D^{i+1} \times S^{n-i-1}~.$$
The trace of the surgery is the cobordism $(W;M,M')$ given by
attaching a $(i+1)$-handle at $S^i \times D^{n-i} \subset M$
$$W~=~M \times I \cup D^{i+1} \times D^{n-i}~.$$
The trace of the surgery on $D^{i+1} \times S^{n-i-1} \subset M'$
is the cobordism $(W';M',M)$ with 
$$W'~=~-W~=~M'\times I \cup D^{i+1} \times D^{n-i}~.$$
In fact, every cobordism of manifolds is a union of the traces of
surgeries.

In terms of homotopy theory the trace $W$ is obtained from $M$ by
attaching an $(i+1)$-cell, and $M'$ is then obtained from $W$ by
detaching an $(n-i)$-cell, with homotopy equivalences
$$W~\simeq~M \cup_x D^{i+1}~\simeq~M' \cup_{x'}D^{n-i}~,$$
with $x:S^i \to M$ the inclusion $S^i\times \{0\} \subset S^i \times D^{n-i}
\subset M$, and similarly for $x':S^{n-i-1} \to M'$.
The immediate homology effect of the surgery is to kill $x \in H_i(M)$, 
$$H_i(W)~=~H_i(M)/\langle x \rangle$$
with $\langle x\rangle \subseteq H_i(M)$ the subgroup generated by $x$.
On the chain level
\begin{itemize}
\item[(i)]  $C(W)$ is chain equivalent to
the algebraic mapping cone ${\mathcal C}(x)$ of a chain map
$x:S^i\Z \to C(M)$ representing $x \in H_i(M)$, where
$$S^i\Z~:~\dots \to 0 \to \Z \to 0 \to \dots~~(\hbox{\rm concentrated
in degree $i$})~,$$
and similarly for $C(W)\simeq {\mathcal C}(x':S^{n-i-1}\Z \to C(M'))$,
\item[(ii)] there is defined a commutative braid of chain homotopy exact
sequences of chain complexes
\vskip-3mm

$$\hskip10pt\xymatrix@!C@C-52pt{
S^i\Z \ar[dr]\ar@/^2pc/[rr]^-{\displaystyle{x}}&&
C(M) \ar[dr] \ar@/^2pc/[rr]^-{\displaystyle{x^*}}&&S^{n-i}\Z \\&
C(W,M \cup M')_{*+1}\ar[ur] \ar[dr] && C(W)
\ar[ur]^-{\displaystyle{y}} \ar[dr]_-{\displaystyle{y'}}&&\\
S^{n-i-1}\Z
\ar[ur]\ar@/_2pc/[rr]_-{\displaystyle{x'}}&&C(M')
\ar[ur]\ar@/_2pc/[rr]_-{\displaystyle{{x'}^*}}&&S^{i+1}\Z
}$$
\vskip2mm

\noindent with $x^*:C(M) \to S^{n-i}\Z$ a chain map
representing the Poincar\'e dual $x^* \in H^{n-i}(M)$ of $x \in H_i(M)$,
and similarly for ${x'}^*$,
\item[(iii)] $C(M')$ is chain equivalent to the dimension shifted
algebraic mapping cone ${\mathcal C}(y)_{*+1}$ of a chain map $y:C(W) \to S^{n-i}\Z$ 
representing a cohomology class $y \in H^{n-i}(W)$ with image the
Poincar\'e dual $x^* \in H^{n-i}(M)$ of $x \in H_i(M)$, and similarly
for $C(M)$.
\end{itemize}
Algebraic surgery gives a precise algebraic model for a chain complex
in the chain homotopy type of $C(M')$, which is obtained from $C(M)$ by
attaching $x$ and detaching $y$.
\medskip

The homology groups $H_*(M),H_*(M'),H_*(W)$ are related by the long
exact sequences
$$\begin{array}{l}
\dots \to H_r(M) \to H_r(W) \to H_r(W,M) \to H_{r-1}(M) \to \dots~,\\[1ex]
\dots \to H_r(M') \to H_r(W) \to H_r(W,M') \to H_{r-1}(M') \to \dots~.
\end{array}$$
It now follows from the excision isomorphisms
$$\begin{array}{l}
H_r(W,M)~=~H_r(D^{i+1},S^i)~=~
\begin{cases}
~\Z&{\rm for}~r=i+1\\
~0&{\rm otherwise}~,
\end{cases}\\[3ex]
H_r(W,M')~=~H_r(D^{n-i},S^{n-i-1})~=~
\begin{cases}
~\Z&{\rm for}~r=n-i\\
~0&{\rm otherwise}
\end{cases}
\end{array}$$
that
$$H_r(M)~=~H_r(W)~=~H_r(M')~{\rm for}~r \neq i,i+1,n-i-1,n-i~.$$
The relationship between $H_r(M),H_r(M'),H_r(W)$ for $r=i,i+1,n-i-1,n-i$
is more complicated, especially in the middle dimensional cases
$n=2i,2i+1$. 

Here are some of the advantages of chain complexes over homology
in describing the algebraic effects of surgery on manifolds.
The chain complex method~:
\begin{itemize}
\item[($\bullet$)] makes it easier to follow the passage from the 
embedding $S^i \times D^{n-i} \subset M$ to the homology $H_*(M')$ 
on the chain level;
\item[($\bullet$)] provides a uniform description for all $n,i$\,;
\item[($\bullet$)] avoids the indeterminacies inherent in exact sequences;
\item[($\bullet$)] works just as well in the non-simply connected case;
\item[($\bullet$)] keeps track of the effect of successive surgeries.
\end{itemize}

Surgery on manifolds is described algebraically by surgery on
chain complexes with symmetric Poincar\'e duality.
The applications of surgery to the classification of manifolds
involve a normal map $(f,b):M \to X$, and only surgeries 
with an extension of $(f,b)$ to a normal map on the trace
$$((g,c);(f,b),(f',b'))~:~(W;M,M') \to X \times ([0,1];\{0\},\{1\})$$
are considered. Surgery on normal maps is described 
algebraically by surgery on chain complexes with quadratic
Poincar\'e duality. The quadratic refinement corresponds to the
additional information carried by the bundle map $b:\nu_M \to \nu_X$.
The formulae for algebraic surgery on symmetric Poincar\'e complexes
are entirely analogous to the formulae for quadratic Poincar\'e complexes.

\subsection{The Principle of Algebraic Surgery}

In its simplest form, the Principle states that for a cobordism of
$n$-dimensional manifolds $(W;M,M')$ the chain homotopy type of
$C(M')$ and the Poincar\'e duality chain equivalence
$$[M'] \cap -~:~C(M')^{n-*}~ \simeq~ C(M')$$ 
can be obtained from 
\begin{itemize}
\item[(i)] the chain homotopy type of $C(M)$, 
\item[(ii)] the Poincar\'e duality chain equivalence 
$$\phi_0~=~[M]\cap -~:~C(M)^{n-*}~ \simeq~ C(M)$$
and the chain homotopy
$$\phi_1~:~(\phi_0)^*~\simeq~\phi_0~:~C(M)^{n-*} \to C(M)$$
determined up to higher chain homotopies by topology,
\item[(iii)] the chain homotopy class of the chain map $j:C(M) \to C(W,M')$ 
induced by the inclusion $M \subset W$,
\item[(iv)] the chain homotopy
$$\delta\phi_0~:~j\phi_0j^*~\simeq~0~:~C(W,M')^{n-*} \to C(W,M')$$
determined up to higher chain homotopies by topology.
\end{itemize} 
The chain complex $C(M')$ is chain equivalent to the chain complex $C'$ 
obtained from $C(M)$ by algebraic surgery, with
$$C'_r~=~C_r(M) \oplus C_{r+1}(W,M') \oplus C^{n-r-1}(W,M')~.$$
See \S3 for formulae for the differentials and Poincar\'e duality of $C'$.

In particular, if $(W;M,M')$ is the trace of a surgery on
$S^i\times D^{n-i} \subset M$ as in \S\ref{vs} then $C(W,M')$ 
is chain equivalent to $S^{n-i}\Z$, and replacing $C(W,M')$ by 
$S^{n-i}\Z$ in the formula for $C'_r$ gives a smaller chain complex 
(also denoted by $C'$) 
$$\begin{array}{l}
C'~:~C_n(M) \to \dots \to C_{n-i}(M) \xrightarrow[]{d \oplus y} C_{n-i-1}(M) \oplus \Z
\xrightarrow[]{d\oplus 0} C_{n-i-2}(M) \to \dots\\[1ex]
\hskip100pt \to C_{i+2}(M) \xrightarrow[]{d \oplus 0} C_{i+1}(M) \oplus \Z 
\xrightarrow[]{d \oplus x} C_i(M) \to \dots \to C_0(M)
\end{array}$$
chain equivalent to $C(M')$. The attaching chain map
$x:S^i \Z \to C(M)$ and the chain map $j:C(M) \to C(W,M') \simeq
S^{n-i}\Z$ in (iii) are determined by the homotopy class of the core
embedding $S^i \times \{0\} \subset M$. 
The detaching chain map $y:{\mathcal C}(x) \to S^{n-i}\Z$ and the
chain homotopy $\delta\phi_0$ in (iv) are determined
by the framing of the core, and are much more subtle.
(See the Examples below).
In this case the algebraic surgery kills the homology class
$x \in H_i(M)$. In the general algebraic context surgery kills entire 
subcomplexes rather than just individual homology classes.

\noindent{\it Example.} The effect of surgery on
$S^0\times D^1\subset M=S^1$ is a double cover of $S^1$.
There are two possibilities:\\
(i) If the two paths $S^0\times D^1 \subset S^1$ move in opposite senses the
effect of the surgery is the trivial double cover 
$M'=S^1 \cup S^1$ of $S^1$, and the trace $(W;M,M')$ is given by
the orientable 
$$W~=~{\rm cl.}(S^2\backslash (D^2 \cup D^2 \cup D^2))~.$$
(ii) If the two paths $S^0 \times D^1 \subset S^1$ move in the same sense the
effect of the surgery is the nontrivial double cover $M''=S^1$
of $S^1$, and the trace $(W';M,M'')$ is given by the nonorientable
$$W'~=~{\rm cl.}(\hbox{\rm M\"obius band}\backslash D^2)~.\eqno{\square}$$

More generally:

\noindent{\it Example.} As usual, let $O(j)$ be the orthogonal group
of $\R^j$. For any map $\omega:S^i \to O(j)$ write $n=i+j$, and
define an embedding
$$e_{\omega}~:~S^i \times D^j \to S^n=S^i\times D^j \cup
D^{i+1} \times S^{j-1}~;~(x,y) \mapsto (x,\omega(x)(y))~.$$
Surgery on $M=S^n$ killing $e_{\omega}$ has effect the $(j-1)$-sphere bundle 
over $S^{i+1}=D^{i+1}\cup D^{i+1}$
$$M'~=~S(\omega)~=~D^{i+1} \times S^{j-1} \cup_{\omega} D^{i+1} \times S^{j-1}$$
of the $j$-plane vector bundle over $S^{i+1}$ 
$$E(\omega)~=~D^{i+1} \times \R^j \cup_{\omega} D^{i+1} \times \R^j$$
classified by $\omega \in \pi_i(O(j))=\pi_{i+1}(BO(j))$.
The trace of the surgery is 
$$W~=~{\rm cl.}(D(\omega)\backslash D^{n+1})~,$$
with 
$$D(\omega)~=~D^{i+1} \times D^j \cup_{\omega} D^{i+1} \times D^j$$
the $j$-disk bundle of $\omega$, which fits into a fibre bundle
$$(D^j,S^{j-1}) \to (D(\omega),S(\omega)) \to S^{i+1}~.$$
Exercise: work out the algebraic effect of the surgery!\hfill\qed

\section{Forms and formations}

The quadratic $L$-groups $L_*(A)$ were originally defined by Wall, with
$L_{2i}(A)$ a Witt-type group of stable isomorphism classes of
nonsingular $(-)^i$-quadratic forms over a ring with involution $A$,
and $L_{2i+1}(A)$ a Whitehead-type group of automorphisms of
$(-)^i$-quadratic forms over $A$ (now replaced by formations).  The
surgery obstruction of a normal map $(f,b):M \to X$ from an
$n$-dimensional manifold $M$ to an $n$-dimensional Poincar\'e complex $X$
$$\sigma_*(f,b) \in L_n(\Z[\pi_1(X)])$$
was defined by first making $(f,b)$ $i$-connected for $n=2i$ (resp. $2i+1$)
by surgery below the middle dimension. The surgery obstruction
is such that $\sigma_*(f,b)=0$ if (and for $n>4$ only if) $(f,b)$ 
is normal bordant to a homotopy equivalence.

Let $A$ be an associative ring with 1, and with an involution
$A \to A;a \mapsto \overline{a}$ satisfying
$$\overline{a+b}~=~\overline{a}+\overline{b}~,~
\overline{ab}~=~\overline{b}\overline{a}~,~
\overline{\overline{a}}~=~a~,~\overline{1}~=~1~.$$
In the applications to topology $A=\Z[\pi]$ is a group ring with the
involution
$$\Z[\pi] \to \Z[\pi]~;~x~=~\sum\limits_{g \in \pi}n_gg
\mapsto \overline{x}~=~\sum\limits_{g \in \pi}n_gg^{-1}~.$$

The dual of a left $A$-module is the left $A$-module 
$$K^*~=~\hbox{\rm Hom}_A(K,A)~,~A \times K^* \to K^*~;~(a,f) \mapsto (x
\mapsto f(x)\overline{a})~.$$
The dual of an $A$-module morphism $f:K \to L$ is the $A$-module morphism
$$f^*~:~L^* \to K^*~;~g \mapsto (x \mapsto g(f(x)))~.$$
For f.g. free $K,L$ identify
$$f^{**}~=~f~:~K^{**}~=~K \to L^{**}~=~L~,$$
using the isomorphism $K \to K^{**};x \mapsto (f \mapsto \overline{f(x)})$
to identify $K=K^{**}$, and similarly for $L$. 

A {\it $(-)^i$-quadratic form} $(K,\lambda,\mu)$ is a f.g.
free $A$-module $K$ together with a $(-)^i$-symmetric form
$$\lambda~=~(-)^i\lambda^*~:~K \to K^*$$
and a function
$$\mu~:~K \to Q_{(-)^i}(A)~=~A/\{a-(-)^i\overline{a}\,\vert\,a \in A\}$$
such that 
$$\lambda(x)(x)~=~\mu(x)+(-)^i\overline{\mu(x)}~,~
\mu(ax)~=~a\mu(x)\overline{a}~,~
\mu(x+y)~=~\mu(x)+\mu(y)+\lambda(x,y)~.$$
The form is {\it nonsingular} if $\lambda:K \to K^*$ is an isomorphism. 

A {\it lagrangian} for a nonsingular $(-)^i$-quadratic form
$(K,\lambda,\mu)$ is a f.g.  free direct summand $L\subset K$ such that
$\lambda(L)(L)=0$, $\mu(L)=0$, and $L=L^{\perp}$, where
$$L^{\perp}~=~\{x \in K\,\vert\,\lambda(x)(L)=0\}~.$$
A nonsingular form admits a lagrangian if and only if it is isomorphic
to the hyperbolic form
$$H_{(-)^i}(L)~=~(L \oplus L^*,\begin{pmatrix} 0 & 1 \cr (-)^i & 0
\end{pmatrix},\mu)$$
with $\mu(x,f)=f(x)$. 

The {\it $2i$-dimensional quadratic $L$-group} $L_{2i}(A)$ is the
Witt group of stable isomorphism classes of nonsingular $(-)^i$-quadratic
forms $(K,\lambda,\mu)$ over $A$, where stability is with respect to
the hyperbolic forms. 

As in Chapter 5 of Wall \cite{Wall} the surgery obstruction of an $i$-connected
$2i$-dimensional normal map $(f,b):M \to X$ is the Witt class
$$\sigma_*(f,b)~=~(K_i(M),\lambda,\mu) \in L_{2i}(\Z[\pi_1(X)])$$
of the kernel nonsingular $(-)^i$-quadratic form $(K_i(M),\lambda,\mu)$
over $\Z[\pi_1(X)]$,
with $\lambda,\mu$ defined by geometric intersection numbers.
An algebraic surgery on $(f,b)$ removing $S^j \times D^{2i-j} \subset M$ 
for $j=i-1$ (resp. $i$) correspond to the algebraic surgery of the
addition (resp.  subtraction) of the hyperbolic form
$H_{(-)^i}(\Z[\pi_1(X)])$ to (resp.  from) the kernel form.

A {\it nonsingular $(-)^i$-quadratic formation} $(K,\lambda,\mu;F,G)$
is a nonsingular $(-)^i$-quadratic form $(K,\lambda,\mu)$ together with
an ordered pair of lagrangians $F,G$. 

The {\it $(2i+1)$-dimensional quadratic $L$-group} $L_{2i+1}(A)$ is the
group of stable isomorphism classes of nonsingular $(-)^i$-quadratic
formations $(K,\lambda,\mu;F,G)$ over $A$, where stability is with respect to
the formations such that either $F,G$ are direct complements in $K$ or
share a common lagrangian complement in $K$.

The surgery obstruction of an $i$-connected
$(2i+1)$-dimensional normal map $(f,b):M \to X$ 
$$\sigma_*(f,b)~=~(K,\lambda,\mu;F,G) \in L_{2i+1}(\Z[\pi_1(X)])$$
is the Witt-type equivalence class of a kernel
$(-)^i$-quadratic formation over $\Z[\pi_1(X)]$ with
$$F\cap G~=~K_{i+1}(M)~~,~~K/(F+G)~=~K_i(M)~.$$
As in Chapter 6 of Wall \cite{Wall} such a kernel formation 
$(K,\lambda,\mu;F,G)$ is obtained by
realizing any finite set $\{x_1,x_2,\dots,x_k\}\subset K_i(M)$ of
$\Z[\pi_1(X)]$-module generators by a high-dimensional Heegaard-type
decomposition of $(f,b)$ as a union of normal maps
$$(f,b)~=~(f_0,b_0) \cup (g,c)~:~M~=~M_0 \cup U \to X~=~X_0 \cup D^{2i+1}$$
with 
$$\begin{array}{l}
(g,c)~:~ (U,\partial U)~=~(\#_k S^i\times D^{i+1},\#_k S^i\times S^i)
\to (D^{2i+1},S^{2i})~,\\[1ex]
F~=~{\rm im}(K_{i+1}(U,\partial U) \to K_i(\partial U))~=~
\Z[\pi_1(X)]^k~,\\[1ex]
G~=~{\rm im}(K_{i+1}(M_0,\partial U) \to K_i(\partial U))~\cong~
\Z[\pi_1(X)]^k~,\\[1ex]
K~=~K_i(\partial U)~=~F \oplus F^*~,~(\lambda,\mu)~=~ \hbox{\rm hyperbolic
$(-)^i$-quadratic form}~.
\end{array}$$

\section{Surgery on symmetric Poincar\'e complexes}

Symmetric Poincar\'e complexes are chain complexes with the Poincar\'e
duality properties of manifolds.  A manifold $M$ determines a symmetric
Poincar\'e complex $(C,\phi)$, such that a surgery on $M$ determines an
algebraic surgery on $(C,\phi)$.  However, not every algebraic surgery
on $(C,\phi)$ can be realized by a surgery on $M$.

Given a f.g. free $A$-module chain complex
$$C~:~ \dots \to C_{r+1} \xrightarrow[]{d} C_r \xrightarrow[]{d} C_{r-1} \to \dots$$
write the dual f.g. free $A$-modules as
$$C^r~=~(C_r)^*~.$$
For any $n\geqslant 0$ let $C^{n-*}$ be the f.g. free $A$-module chain complex with
$$d_{C^{n-*}}~=~(-)^rd_C^*~:~(C^{n-*})_r~=~C^{n-r} \to
(C^{n-*})_{r-1}~=~C^{n-r+1} ~.$$
The duality isomorphisms 
$$T~:~{\rm Hom}_A(C^p,C_q) \to {\rm Hom}_A(C^q,C_p)~;~
\phi \mapsto (-)^{pq}\phi^*$$
are involutions with the property that the dual of a chain map
$f:C^{n-*} \to C$ is a chain map $Tf:C^{n-*} \to C$, with $T(Tf)=f$.

The algebraic mapping cone ${\mathcal C}(f)$ of a chain map $f:C \to D$
is the chain complex with
$$d_{C(f)}~=~\begin{pmatrix} d_D & (-1)^rf \cr 0 & d_C \end{pmatrix}~:~
C(f)_r ~=~D_r \oplus C_{r-1} \to C(f)_{r-1} ~=~D_{r-1} \oplus C_{r-2}~.$$
\indent
An {\it $n$-dimensional symmetric complex} $(C,\phi)$ over $A$ is a 
f.g. free $A$-module chain complex
$$C~:~C_n \xrightarrow[]{d_C} C_{n-1} \to \dots \to C_1\xrightarrow[]{d_C} C_0$$
together with a collection of $A$-module morphisms
$$\phi~=~\{\phi_s:C^{n-r+s} \to C_r\,\vert\,s \geqslant 0\}$$
such that
$$\begin{array}{c}
d_C\phi_s+(-1)^r\phi_sd_C^*+(-1)^{n+s-1}(\phi_{s-1}+(-1)^sT\phi_{s-1})~=~0~:~
C^{n-r+s-1} \to C_r\\[1ex]
(s \geqslant 0,\phi_{-1}=0)~.
\end{array}$$
Thus $\phi_0:C^{n-*} \to C$ is a chain map, $\phi_1$ is a chain homotopy
$\phi_1:\phi_0\simeq T\phi_0:C^{n-*} \to C$, and so on $\dots$~. 
More intrinsically, $\phi$ is an $n$-dimensional cycle in the 
$\Z$-module chain complex
$${\rm Hom}_{\Z[\Z_2]}(W,{\rm Hom}_A(C^*,C))$$
with 
$$W~:~\dots \to \Z[\Z_2] \xrightarrow[]{1+T} \Z[\Z_2] \xrightarrow[]{1-T} \Z[\Z_2]
 \xrightarrow[]{1+T} \Z[\Z_2] \xrightarrow[]{1-T} \Z[\Z_2]$$
the free $\Z[\Z_2]$-module resolution of $\Z$.
The symmetric complex $(C,\phi)$ is {\it Poincar\'e}
if the chain map $\phi_0:C^{n-*} \to C$ is a chain equivalence.

\noindent{\it Example.} (Mishchenko \cite{Mishchenko})
An $n$-dimensional manifold $M$ and a regular covering $\widetilde{M}$ 
with group of covering translations $\pi$ determine
an $n$-dimensional symmetric Poincar\'e complex over $\Z[\pi]$
$(C(\widetilde{M}),\phi)$ with
$$\phi_0~=~[M]\cap -~:~C(\widetilde{M})^{n-*} \to C(\widetilde{M})$$
the Poincar\'e duality chain equivalence.
The higher chain homotopies $\phi_1,\phi_2,\dots$ are determined
by an equivariant analogue of the construction of the Steenrod squares.
\hfill\qed

An {\it $(n+1)$-dimensional symmetric pair} 
$(j:C \to D,(\delta\phi,\phi))$ 
is an $n$-dimensional symmetric complex $(C,\phi)$ together with a
chain map $j:C \to D$ to an $(n+1)$-dimensional f.g. free
$A$-module chain complex $D$ and $A$-module morphisms
$\delta\phi=\{\delta\phi_s:D^{n+1-r+s} \to D_r\,\vert\,s \geqslant 0\}$
such that
$$\begin{array}{c}
j\phi_sj^*~=~d_D\delta\phi_s+(-)^r\delta\phi_sd_D^*+
(-)^{n+s+1}(\delta\phi_{s-1}+(-)^sT\delta\phi_{s-1})~:~
D^{n+1-r-s} \to D_r\\[1ex]
(s \geqslant 0,\delta\phi_{-1}=0)~.
\end{array}$$
The pair is {\it Poincar\'e} if the chain map
$$\begin{pmatrix}\delta\phi_0 \\ \phi_0j^* \end{pmatrix}~:~
D^{n+1-*} \to {\mathcal C}(j)$$
is a chain equivalence, in which case $(C,\phi)$ is a
$n$-dimensional symmetric Poincar\'e complex. 

A {\it cobordism} of $n$-dimensional symmetric Poincar\'e complexes
$(C,\phi)$, $(C',\phi')$ is an $(n+1)$-dimensional symmetric Poincar\'e
pair of the type $(C \oplus C' \to D,(\delta\phi,\phi\oplus -\phi'))$. 
Symmetric complexes $(C,\phi)$, $(C',\phi')$ are {\it homotopy
equivalent} if there exists a cobordism with $C \to D$, $C' \to D$
chain equivalences.

\noindent{\it Example.} A 0-dimensional symmetric complex $(C,\phi)$
is a f.g. free $A$-module $C_0$ together with a symmetric form $\phi_0$
on $C^0$. The complex is Poincar\'e if and only if the form is nonsingular.
Two 0-dimensional symmetric Poincar\'e complexes $(C,\phi)$, $(C',\phi')$
are cobordant if and only if the forms $(C^0,\phi_0)$, $({C'}^0,\phi'_0)$
are Witt-equivalent,
i.e. become isomorphic after stabilization with metabolic forms
$(L \oplus L^*, \begin{pmatrix} \lambda & 1 \\ 1 & 0
\end{pmatrix})$.\hfill\qed

\noindent{\it Example.} 
An $(n+1)$-dimensional manifold with boundary $(W,\partial W)$ and
cover $(\widetilde{W},\partial \widetilde{W})$
determines an $(n+1)$-dimensional symmetric Poincar\'e pair 
$(j:C(\partial \widetilde{W}) \to C(\widetilde{W}),(\delta\phi,\phi))$
over $\Z[\pi]$ with
$$\begin{pmatrix} \delta\phi_0 \cr \phi_0j^* \end{pmatrix}~=~[W] \cap -~:~
D^{n+1-*}~=~C(\widetilde{W})^{n-*} \to {\mathcal C}(j)~=~C(\widetilde{W},\partial
\widetilde{W})$$
the Poincar\'e-Lefschetz duality chain equivalence.\hfill\qed

The {\it data} for {\it algebraic surgery} on an $n$-dimensional 
symmetric Poincar\'e complex $(C,\phi)$ is an $(n+1)$-dimensional
symmetric pair $(j:C \to D,(\delta\phi,\phi))$. The {\it effect}
of the algebraic surgery is the $n$-dimensional symmetric Poincar\'e
complex $(C',\phi')$ with
$$\begin{array}{l}
d_{C'}~=~\begin{pmatrix} d_C &0 & (-)^{n+1}\phi_0j^* \\
(-)^rj & d_D & (-)^r \delta\phi_0 \\
0 & 0 & d_D^* \end{pmatrix}~:\\[4ex]
\hskip25pt 
C'_r~=~C_r \oplus D_{r+1} \oplus D^{n-r+1}
\to C'_{r-1}~=~C_{r-1} \oplus D_r \oplus D^{n-r+2}~,\\[2ex]
\phi'_0~=~\begin{pmatrix} \phi_0 &  0 & 0 \\
(-)^{n-r}jT\phi_1  & (-)^{n-r}T\delta\phi_1 & 0 \\
0 & 1 & 0 \end{pmatrix}~:\\[4ex]
\hskip25pt
{C'}^{n-r}~=~C^{n-r} \oplus D^{n-r+1} \oplus D_{r+1} 
\to C'_r~=~C_r \oplus D_{r+1} \oplus D^{n-r+1}~,\\[1ex]
\phi'_s~=~\begin{pmatrix} \phi_s &  0 & 0 \\
(-)^{n-r}jT\phi_{s+1}  & (-)^{n-r+s}T\delta\phi_{s+1} & 0 \\
0 & 0 & 0 \end{pmatrix}~:\\[4ex]
\hskip25pt 
{C'}^{n-r+s}~=~C^{n-r+s} \oplus D^{n-r+s+1} \oplus D_{r-s+1} 
\to C'_r~=~C_r \oplus D_{r+1} \oplus D^{n-r+1}~~(s \geqslant 1)~.
\end{array}$$
\noindent{\it Remark.}
The appearance of the chain homotopy $\phi_1:\phi_0\simeq T\phi_0$ in
the formula for the Poincar\'e duality chain equivalence $\phi'_0$ is a
reason for taking account of $\phi_1$.  The appearance of the higher
chain homotopy $\phi_2:\phi_1 \simeq T\phi_1$ in the formula for
$\phi'_1$ is a reason for taking account of $\phi_2$.  And so on
$\dots$.\hfill\qed
\medskip

The {\it trace} of an algebraic surgery is the $(n+1)$-dimensional
symmetric Poincar\'e cobordism between $(C,\phi)$ and $(C',\phi')$
$$((f~f'):C \oplus C' \to D',(0,\phi \oplus -\phi'))$$
defined by
$$\begin{array}{l}
d_{D'}~=~\begin{pmatrix} d_C & (-)^{n+1}\phi_0j^* \\
0 & d_D^* \end{pmatrix}~:~D'_r~=~C_r  \oplus D^{n-r+1}
\to D'_{r-1}~=~C_{r-1}  \oplus D^{n-r+2}~,\\[2ex]
f~=~\begin{pmatrix} 1  \\ 0  \end{pmatrix}~:~
C_r  \to D'_r~=~C_r  \oplus D^{n-r+1}~,\\[2ex]
f'~=~\begin{pmatrix} 1 & 0 & 0 \\ 0 & 0 & 1 \end{pmatrix}~:~
C'_r~=~C_r \oplus D_{r+1} \oplus D^{n-r+1} \to
D'_r~=~C_r  \oplus D^{n-r+1}~.
\end{array}$$

\noindent{\it Theorem.} (Ranicki \cite{Ranicki1980I}) 
{\it The cobordism of symmetric
Poincar\'e complexes is the equivalence relation generated by homotopy
equivalence and algebraic surgery.}  \\
{\it Proof.} Homotopy equivalent complexes are cobordant, by definition.
Surgery equivalent complexes are cobordant by the trace construction.\\
\indent Conversely, suppose given a cobordism of $n$-dimensional 
symmetric Poincar\'e complexes
$$\Gamma~=~((f~f'):C \oplus C' \to D,(\delta\phi,\phi\oplus -\phi'))~.$$
Let
$$\overline{\Gamma}~=~
((\overline{f}~\overline{f}'):C \oplus \overline{C}' \to 
\overline{D},(0,\phi\oplus -\phi'))$$
be the trace of the algebraic surgery on $(C,\phi)$ with 
data $(j:C \to {\mathcal C}(f'),(\delta\phi/\phi',\phi))$ given by
$$\begin{array}{l}
j~=~\begin{pmatrix} f \\ 0 \end{pmatrix}~:~
C_r \to {\mathcal C}(f')_r~=~D_r \oplus C'_{r-1}~,\\[1ex]
(\delta\phi/\phi')_s~=~
\begin{pmatrix}  \delta\phi_s & (-)^sf'\phi'_s \\
0 & (-)^{n-r+s}T\phi'_{s-1} \end{pmatrix}~:\\[2ex]
\hskip25pt
{\mathcal C}(f')^{n-r+s+1}~=~D^{n-r+s+1} \oplus {C'}^{n-r+s}
\to {\mathcal C}(f')_r ~=~D_r \oplus C'_{r-1}~~(s \geqslant 0,\phi'_{-1}=0)~.
\end{array}$$
The $A$-module morphisms
$$\begin{array}{l}
g~=~(0~~0~~1~~0~~0)~:~\overline{C}'_r~=~
C_r \oplus D_{r+1} \oplus C'_r \oplus D^{n-r+1} \oplus {C'}^{n-r}
\to C'_r~,\\[1ex]
h~=~(f~~\delta\phi_0~f'\phi'_0)~:~\overline{D}_r~=~C_r \oplus D^{n-r+1} \oplus {C'}^{n-r}
\to D_r
\end{array}$$
define a homotopy equivalence 
$(h,1_C \oplus g):\overline{\Gamma} \to \Gamma$.\hfill\qed

\noindent {\it Symmetric Surgery Principle.} 
For any $(n+1)$-dimensional cobordism $(W;M,M')$ 
and regular cover $(\widetilde{W};\widetilde{M},\widetilde{M}')$
with group $\pi$ the symmetric
Poincar\'e complex $(C(\widetilde{M}'),\phi')$ is homotopy equivalent to 
the effect of algebraic surgery on $(C(\widetilde{M}),\phi)$
with data 
$$(j:C(\widetilde{M}) \to C(\widetilde{W},\widetilde{M}'),(\delta\phi',\phi))~.$$
\noindent{\it Proof.} The manifold cobordism determines a cobordism
of $n$-dimensional symmetric Poincar\'e complexes over $\Z[\pi]$
$$\Gamma~=~(C(\widetilde{M}) \oplus C(\widetilde{M}') \to
C(\widetilde{W}),(\delta\phi,\phi\oplus -\phi'))~.$$
Now apply the Theorem to $\Gamma$, with $\delta\phi'=\delta\phi/\phi'$.\hfill$\square$
\medskip

\noindent{\it Example.}
If $(W;M,M')$ is the trace of a surgery on $S^i \times D^{n-i} \subset M$
then 
$$C(\widetilde{W},\widetilde{M}')~ \simeq~ S^{n-i}\Z[\pi]$$
is concentrated in dimension $(n-i)$, and the effect is to kill
the spherical (co)homology class 
$$j~=~[S^i] \in H^{n-i}(\widetilde{M})~\cong~H_i(\widetilde{M})~.$$
The embedding $S^i \subset M$ determines $j: C(\widetilde{M})
\to S^{n-i}\Z[\pi]$, and the choice of extension to an embedding
$S^i \times D^{n-i} \subset M$ determines $\delta\phi'$.
\hfill\qed
\medskip

The cobordism groups $L^n(A)$ ($n
\geqslant 0$) start with the symmetric Witt group $L^0(A)$.  The symmetric
signature map from manifold bordism to symmetric Poincar\'e bordism
$$\sigma^*~:~\Omega_n(X) \to L^n(\Z[\pi_1(X)])~;~M \mapsto
(C(\widetilde{M}),\phi)$$
is a generalization of the signature map
$$\sigma^*~:~\Omega_{4k} \to L^{4k}(\Z)~=~\Z~;~M \mapsto 
{\rm signature}(H^{2k}(M),{\rm intersection~form})~.$$
Although the symmetric signature maps $\sigma^*$ are not isomorphisms 
in general, they do provide many invariants. The symmetric and
quadratic $L$-groups only differ in 8-torsion :
\begin{itemize}
\item[(i)] the symmetrization maps
$$1+T~:~L_n(A) \to L^n(A)~;~(C,\psi) \mapsto (C,(1+T)\psi)$$
are isomorphisms modulo 8-torsion,
\item[(ii)] if $1/2 \in A$ the symmetrization maps are isomorphisms.
\end{itemize}

\section{Surgery on quadratic Poincar\'e complexes}

Quadratic Poincar\'e complexes are chain complexes with the Poincar\'e
duality properties of kernels of normal maps. 
The quadratic Poincar\'e analogues of cobordism and surgery are
defined by analogy with the symmetric case. Although there are
many similarities between the quadratic and symmetric theories,
there is one essential difference : the quadratic Poincar\'e 
cobordism groups are the Wall surgery obstruction groups $L_*(A)$,
so for $A=\Z[\pi]$ every element is geometrically significant.

An {\it $n$-dimensional quadratic complex} $(C,\psi)$ over $A$ is a 
f.g. free $A$-module chain complex $C$ together with a collection
of $A$-module morphisms
$$\psi~=~\{\psi_s:C^{n-r-s} \to C_r\,\vert\,s \geqslant 0\}$$
such that
$$d_C\psi_s+(-1)^r\psi_sd_C^*+(-1)^{n-s-1}(\psi_{s+1}+(-1)^{s+1}T\psi_{s+1})~
=~0~:~C^{n-r-s-1} \to C_r~~(s \geqslant 0)~.$$
More intrinsically, $\psi$ is an $n$-dimensional cycle in the 
$\Z$-module chain complex
$$W \otimes_{\Z[\Z_2]}{\rm Hom}_A(C^*,C)$$
with $W$ the free $\Z[\Z_2]$-module resolution of $\Z$ (as above).
The quadratic complex $(C,\psi)$ is {\it Poincar\'e}
if the chain map $(1+T)\psi_0:C^{n-*} \to C$ is a chain equivalence. 
A quadratic complex $(C,\psi)$ determines the symmetric complex $(C,\phi)$
with $\phi_0=(1+T)\psi_0$, $\phi_s=0$ $(s \geqslant 1)$.

\noindent{\it Example.} (\cite{Ranicki1980II})
An $n$-dimensional normal map $(f,b):M \to X$ and a regular covering
$\widetilde{X}$ of $X$ with group of covering translations $\pi$
determine a {\it kernel} $n$-dimensional quadratic Poincar\'e complex
$(C,\psi)$ over $\Z[\pi]$ with $C={\mathcal C}(f^{\,!})$ the algebraic
mapping cone of the Umkehr chain map
$$f^{\,!}~:~C(\widetilde{X})~\simeq~ C(\widetilde{X})^{n-*}
\xrightarrow[]{f^*} C(\widetilde{M})^{n-*}~ \simeq~ C(\widetilde{M})$$
and $(1+T)\psi_0=[M] \cap -:C^{n-*} \to C$
the Poincar\'e duality chain equivalence.
It follows from $f_*[M]=[X] \in H_m(X)$ ($f$ is degree 1) that there
exists a chain homotopy
$ff^{\,!}\simeq 1:C(\widetilde{X}) \to C(\widetilde{X})$.
The homology $\Z[\pi]$-modules of $C$ are thus the kernels of $f$
$$H_*(C)~=~K_*(M)~=~
\hbox{\rm ker}(f_*:H_*(\widetilde{M}) \to H_*(\widetilde{X}))~,$$
such that 
$$H_*(\widetilde{M})~=~K_*(M) \oplus H_*(\widetilde{X})~.\eqno{\square}$$

An {\it $(n+1)$-dimensional quadratic pair} $(j:C \to D,(\delta\psi,\psi))$ 
is an $n$-dimensional quadratic complex $(C,\psi)$ together with a
chain map $j:C \to D$ to an $(n+1)$-dimensional f.g. free
$A$-module chain complex $D$ and $A$-module morphisms
$$\delta\psi~=~\{\delta\psi_s:D^{n+1-r-s} \to D_r\,\vert\,s \geqslant 0\}$$
such that
$$j\psi_sj^*~=~d_D\delta\psi_s+(-)^r\delta\psi_sd_D^*+
(-)^{n+s+1}(\delta\psi_{s+1}+(-)^sT\delta\psi_{s+1})~:~
D^{n+1-r-s} \to D_r~~(s \geqslant 0)~.$$
The pair is {\it Poincar\'e} if the chain map
$$\begin{pmatrix} (1+T)\delta\psi_0 \cr
(1+T)\psi_0j^*  \end{pmatrix}~:~D^{n+1-*} \to {\mathcal C}(j)$$
is a chain equivalence, in which case $(C,\psi)$ is a
$n$-dimensional quadratic Poincar\'e complex. 

A {\it cobordism} of $n$-dimensional quadratic Poincar\'e
complexes $(C,\psi)$, $(C',\psi')$ is an $(n+1)$-dimensional quadratic
Poincar\'e pair of the type $(C \oplus C' \to D,(\delta\psi,\psi\oplus
-\psi'))$.  Quadratic complexes $(C,\psi)$, $(C',\psi')$ are
{\it homotopy equivalent} if there exists a cobordism with $C \to
D$, $C' \to D$ chain equivalences.

\noindent{\it Example.} 
An $(n+1)$-dimensional normal map of pairs $(g,c):(W,\partial W) \to
(Y,\partial Y)$ determines a kernel $(n+1)$-dimensional quadratic Poincar\'e
pair over $\Z[\pi]$ $(j:C(\partial g^{\,!}) \to C(g^{\,!}),(\delta\psi,\psi))$ with 
$$\begin{pmatrix}
(1+T)\delta\psi_0j^* \cr (1+T)\psi_0 \end{pmatrix}~=~[W]\cap -~:~
C(g^{\,!})^{n+1-*} \to {\mathcal C}(j)$$ 
the Poincar\'e-Lefschetz duality chain equivalence.\hfill\qed

The {\it data} for {\it algebraic surgery} on an $n$-dimensional 
quadratic Poincar\'e complex $(C,\psi)$ is an $(n+1)$-dimensional
quadratic pair $(j:C \to D,(\delta\psi,\psi))$. The {\it effect}
of the algebraic surgery is the $n$-dimensional quadratic Poincar\'e
complex $(C',\psi')$ with
$$\begin{array}{l}
d_{C'}~=~\begin{pmatrix} d_C &0 & (-)^{n+1}(1+T)\psi_0j^* \\
(-)^rj & d_D & (-)^r (1+T)\delta\psi_0 \\
0 & 0 & d_D^* \end{pmatrix}~:\\[4ex]
\hskip25pt 
C'_r~=~C_r \oplus D_{r+1} \oplus D^{n-r+1}
\to C'_{r-1}~=~C_{r-1} \oplus D_r \oplus D^{n-r+2}~,\\[2ex]
\psi'_0~=~\begin{pmatrix} \psi_0 &  0 & 0 \\
0 & 0 & 0 \\
0 & 1 & 0 \end{pmatrix}~:\\[4ex]
\hskip25pt
{C'}^{n-r}~=~C^{n-r} \oplus D^{n-r+1} \oplus D_{r+1} 
\to C'_r~=~C_r \oplus D_{r+1} \oplus D^{n-r+1}~,\\[2ex]
\psi'_s~=~\begin{pmatrix} \psi_s & (-)^{r+s}T\psi_{s-1}j^* & 0 \\
0  & (-)^{n-r-s+1}T\delta\psi_{s-1} & 0 \\
0 & 0 & 0 \end{pmatrix}~:\\[4ex]
\hskip25pt 
{C'}^{n-r-s}~=~C^{n-r-s} \oplus D^{n-r-s+1} \oplus D_{r+s+1} 
\to C'_r~=~C_r \oplus D_{r+1} \oplus D^{n-r+1}~~(s \geqslant 1)~.
\end{array}$$
The trace of the algebraic surgery is an $(n+1)$-dimensional
quadratic Poincar\'e cobordism
$((f~f'):C \oplus C' \to D',(\delta\psi',\psi \oplus -\psi'))$
defined by analogy with the symmetric case. As in the symmetric case :

\noindent{\it Theorem.} (Ranicki \cite{Ranicki1980I}) {\it
The cobordism of quadratic
Poincar\'e complexes is the equivalence relation generated by homotopy
equivalence and algebraic surgery.} \hfill\qed

\noindent {\it Quadratic Surgery Principle.} 
For a bordism of $n$-dimensional normal maps
$$((g,c);(f,b),(f',b'))~:~(W;M,M') \to X \times ([0,1];\{0\},\{1\})$$
the quadratic Poincar\'e complex $(C({f'}^{\,!}),\psi')$ is homotopy 
equivalent to the effect of algebraic surgery on $(C(f^{\,!}),\psi)$
with data 
$(C(f^{\,!}) \to C(g^{\,!},{f'}^{\,!}),(\delta\psi,\psi))$.
\hfill\qed

\noindent{\it Example.} If $((g,c);(f,b),(f',b'))$
is the trace of a surgery on $S^i \times D^{n-i} \subset M$ then 
$$C(g^{\,!},{f'}^{\,!})~ \simeq~ S^{n-i}\Z[\pi]$$
is concentrated in dimension $(n-i)$.\hfill\qed

A $n$-dimensional quadratic Poincar\'e complex $(C,\psi)$ is
{\it highly-connected} if it is homotopy equivalent to
a complex (also denoted $(C,\psi)$) with
$$\begin{array}{ll}
C~:~\dots \to 0 \to C_i \to 0 \to \dots &\hbox{if $n=2i$} \\[1ex]
C~:~\dots \to 0 \to C_{i+1} \to C_i \to 0 \to \dots&\hbox{if $n=2i+1$.} 
\end{array}$$

\noindent{\it Example.} 
(i) The quadratic kernel $(C,\psi)$ of
an $n$-dimensional normal map $(f,b):M \to X$ is highly-connected
if and only if $f:M \to X$ is $i$-connected, that is $\pi_r(f)=0$
for $r \leqslant i$.\\
(ii) The quadratic Poincar\'e kernel $(C,\psi)$ of an $i$-connected
$2i$-dimensional normal map $(f,b):M \to X$ is essentially the same as the
geometric $(-)^i$-quadratic intersection form $(K_i(M),\lambda,\mu)$
of Wall \cite{Wall}, with 
$$\begin{array}{l}
\lambda~=~(1+T)\psi_0~:~H^i(C)~=~K^i(M) \to H_i(C)~=~K_i(M)~\cong~H^i(C)^*~,\\[1ex]
\mu(x)~=~\psi_0(x)(x) \in Q_{(-)^i}(\Z[\pi_1(X)])~.
\end{array}
$$
For $i \geqslant 3$ an element $x \in K_i(M)$ can be killed by a geometric surgery 
if and only if $\mu(x)=0$, if and only if there exists algebraic surgery
data $(x:C \to S^i\Z[\pi_1(X)],(\delta\psi,\psi))$.
The effect of the surgery is a normal map $(f',b'):M' \to X$ with
quadratic Poincar\'e kernel $(C',\psi')$ such that
$$C'~:~\dots \to 0 \to \Z[\pi_1(X)] \xrightarrow[]{x} K_i(M)
\xrightarrow[]{x^*\lambda} \Z[\pi_1(X)] \to 0 \to \dots~.\eqno{\square}$$

\noindent{\it Theorem.} (Ranicki \cite{Ranicki1980I})\\
(i) {\it Every $n$-dimensional quadratic Poincar\'e complex $(C,\psi)$ is
cobordant to a highly-connected complex.}\\
(ii) {\it The cobordism group of $n$-dimensional quadratic complexes over $A$
is isomorphic to $L_n(A)$, with the 4-periodicity isomorphisms given by}
$$L_n(A) \to L_{n+4}(A)~;~(C,\psi) \mapsto (C_{*-2},\psi)~.$$
\noindent{\it Proof}: (i) Let $n=2i$ or $2i+1$. 
Let $D$ be the quotient complex of $C$ with $D_r=C_r$ for $r>n-i$,
and let $j:C \to D$ be the projection. The effect of algebraic
surgery on $(C,\psi)$ with data $(j:C \to D,(0,\psi))$ 
is homotopy equivalent to a highly-connected complex $(C',\psi')$. \\
(ii) ($n=2i$) A highly-connected $2i$-dimensional quadratic Poincar\'e
complex $(C,\psi)$ is essentially the same as a nonsingular $(-)^i$-quadratic
form $(C^0,\psi_0)$. The relative version of (i) shows that 
a null-cobordism of $(C,\psi)$ is essentially the same as 
an isomorphism of forms 
$$(C^0,\psi_0) \oplus {\rm hyperbolic}~\cong~{\rm hyperbolic}~,$$
which is precisely the condition for $(C^0,\psi_0)=0 \in L_{2i}(A)$.\\
(ii) ($n=2i+1$) A highly-connected $(2i+1)$-dimensional quadratic Poincar\'e
complex $(C,\psi)$ is essentially the same as a nonsingular $(-)^i$-quadratic
formation. See \cite{Ranicki1980I} for further details.\\
\hbox to\hsize{\hfill\qed}

\noindent{\it Instant surgery obstruction for $n=2i$.}
A $2i$-dimensional quadratic Poincar\'e
complex $(C,\psi)$ is cobordant to the highly-connected complex
$(C',\psi')$ with
$$({C'}^i,\psi'_0)~=~\bigg(\hbox{\rm coker}\big(
\begin{pmatrix}
d^* & 0 \\ 
(-)^{i+1}(1+T)\psi_0 & d 
\end{pmatrix}:
C^{i-1} \oplus C_{i+2} \to C^i \oplus C_{i+1}\big),
\begin{pmatrix} \psi_0 & d \\ 0 & 0 \end{pmatrix}\bigg)~.$$
Thus if $(C,\psi)$ is the quadratic Poincar\'e kernel of a
$2i$-dimensional normal map $(f,b):M \to X$ then $({C'}^i,\psi'_0)$ is
a nonsingular $(-)^i$-quadratic form representing the surgery
obstruction $\sigma_*(f,b) \in L_{2i}(\Z[\pi_1(X)])$ (without
preliminary geometric surgeries below the middle dimension).\hfill\qed

See \S I.4 of \cite{Ranicki1980I} for the instant surgery
obstruction formation in the case $n=2i+1$.

\section{The localization exact sequence}

For any morphism of rings with involution $f:A \to B$ there is defined
an exact sequence of $L$-groups
$$\dots \to L_n(A) \xrightarrow[]{f} L_n(B) \to L_n(f) \to L_{n-1}(A)
\to \dots$$
with the relative $L$-group $L_n(f)$ the cobordism groups of pairs
$$\begin{array}{l}
\hbox{($(n-1)$-dimensional quadratic Poincar\'e complex $(C,\psi)$
over $A$,}\\[1ex]
\hskip50pt \hbox{$n$-dimensional quadratic Poincar\'e pair
$(B\otimes_AC \to D,(\delta\psi,1\otimes \psi))$ over $B$)~.}
\end{array}$$
Algebraic surgery provides a particularly useful expression for
the relative $L$-groups $L_*(A \to S^{-1}A)$ of the localization map
$A \to S^{-1}A$ inverting a multiplicatively closed subset $S \subset A$
of central non-zero divisors. 

\noindent{\it Localization exact sequence.} (Ranicki \cite{Ranicki1982})
The relative $L$-group $L_n(A \to S^{-1}A)$ is isomorphic to the
cobordism group $L_n(A,S)$ of $(n-1)$-dimensional quadratic Poincar\'e
complexes over $A$ which are $S^{-1}A$-acyclic.\\
\noindent{\it Proof} Clearing denominators it is possible
to lift every quadratic Poincar\'e pair over $S^{-1}A$ 
as above to an $n$-dimensional quadratic pair 
$(C \to D',(\delta\psi',\psi))$ over $A$. This is the data for
algebraic surgery on $(C,\psi)$ with effect a cobordant
$(n-1)$-dimensional quadratic Poincar\'e complex $(C',\psi')$
over $A$ which is $S^{-1}A$-acyclic (i.e. $H_*(S^{-1}A\otimes_AC')=0$).
\hfill\qed

The localization exact sequence 
$$\dots \to L_n(A) \to L_n(S^{-1}A) \to L_n(A,S) \to L_{n-1}(A) \to \dots$$
and its extensions to noncommutative
localization and to symmetric $L$-theory have many applications to the
computation of $L$-groups, as well as to surgery on submanifolds (cf. 
Ranicki \cite{Ranicki1998}).

\noindent{\it Example.}  The localization of $A=\Z$ inverting
$S=\Z\backslash \{0\} \subset A$ is $S^{-1}A=\Q$.
See Chapter 4 of \cite{Ranicki1982} for detailed accounts of 
the way in which the classification of quadratic forms
over $\mathbb Q$ is combined with the localization exact sequence
$$\begin{array}{l}
\dots \to L^n(\Z) \to L^n(\Q) \to L^n(\Z,S) \to L^{n-1}(\Z) \to \dots~,\\[1ex]
\dots \to L_n(\Z) \to L_n(\Q) \to L_n(\Z,S) \to L_{n-1}(\Z) \to \dots
\end{array}$$
to give 
$$\begin{array}{l}
L^n(\Z)~=~\begin{cases} 
\Z&\text{(signature)}\\ 
\Z_2&\text{(deRham invariant)}\\ 
0&\\ 
0&
\end{cases}
\text{if}~n \equiv \begin{cases} 
0\\ 1\\ 2\\ 3
\end{cases}
(\bmod\,4)~,\\[6ex]
L_n(\Z)~=~\begin{cases} 
\Z&\text{(signature)/8}\\ 0&\\ 
\Z_2&\text{(Arf invariant)}\\ 
0&
\end{cases}\hskip26pt
\text{if}~n \equiv \begin{cases} 
0\\ 1\\ 2\\ 3
\end{cases}
(\bmod\,4)~.
\end{array}$$
$$\eqno{\square}$$


\newpage
\addcontentsline{toc}{section}{References}
\providecommand{\bysame}{\leavevmode\hbox to3em{\hrulefill}\thinspace}


\begin{thebibliography}{99}
\bibitem{Browder}
      W. Browder, {\it Surgery on simply-connected manifolds},
      Ergebnisse der Mathematik 65, Springer (1972)
\bibitem{KervaireMilnor}
      M. Kervaire and J. Milnor, {\it Groups of homotopy spheres I.},
      Ann. of Maths. 77, 504--537 (1963)
\bibitem{Mishchenko}
      A.S.  Mishchenko, {\it Homotopy invariants of non--simply
      connected manifolds, III.  Higher signatures}, Izv.  Akad.  Nauk SSSR,
      ser.  mat.  35, 1316--1355 (1971)
\bibitem{Ranicki1980I}
      A.A. Ranicki, {\it The algebraic theory of surgery I. Foundations},
      Proc. Lond. Math. Soc. (3) 40, 87--192 (1980)
\bibitem{Ranicki1980II}
      \bysame, {\it The algebraic theory of surgery II. Applications
      to topology},
      Proc. Lond. Math. Soc. (3) 40, 193--283 (1980)
\bibitem{Ranicki1982}
      \bysame, {\it Exact sequences in the algebraic theory of surgery},
      Mathematical Notes 26, Princeton (1982) 
      http://www.maths.ed.ac.uk/$\widetilde{~}$aar/books/exact.pdf
\bibitem{Ranicki1998}
      \bysame, {\it High dimensional knot theory}, Springer Monograph,
      Springer (1998)
\bibitem{Ranicki2001}
      \bysame, {\it An introduction to algebraic surgery},
      in Surveys on Surgery Theory, Volume 2, Annals of Mathematics Studies
      149, Princeton, 81--164 (2001)\hfil\break
      e-print http://arXiv.org/abs/math.AT/0008071
\bibitem{Ranicki2002}
      \bysame, {\it Algebraic Poincar\'e cobordism},
Proc. 1999 Conference for 60th birthday of J. Milgram , Stanford,
Contemporary Mathematics 279, A.M.S., 213--255 (2001)\hfil\break
e-print http://arXiv.org/abs/math.AT/0008228
\bibitem{Wall}
      C.T.C. Wall, {\it Surgery on compact manifolds},
      1st Edition, Academic Press (1970),
      2nd Edition, Mathematical Surveys and Monographs 69, A.M.S. (1999)
\end{thebibliography}
\end{document}